\documentclass[graybox]{svmult}

\usepackage{mathtools}

\usepackage{type1cm}        
\usepackage{makeidx}         
\usepackage{graphicx}        
\usepackage{multicol}        
\usepackage[bottom]{footmisc}

\usepackage{newtxtext}       %
\usepackage[varvw]{newtxmath}       

\usepackage{cleveref}
\usepackage{wrapfig}
\usepackage{tikz}
\usepackage{pgfplots}
\pgfplotsset{compat=1.18}
\usepackage{tikz}
\usepackage{booktabs}
\usepackage{todonotes}
\usepackage{siunitx}
\makeindex             

\begin{document}
\title*{Domain decomposition architectures and Gauss--Newton training for physics-informed neural networks}
\titlerunning{Domain decomposition and Gauss--Newton training for PINNs}
\author{Alexander Heinlein\orcidID{0000-0003-1578-8104} \and Taniya Kapoor\orcidID{0000-0002-6361-446X}}
\institute{Alexander Heinlein \at Delft Institute of Applied Mathematics, Delft University of Technology, Mekelweg 4, Delft, 2628 CD, the Netherlands, \email{a.heinlein@tudelft.nl}
  \and Taniya Kapoor \at Artificial Intelligence Group, Wageningen University \& Research, Wageningen, The Netherlands, \email{taniya.kapoor@wur.nl}}
%
%
\maketitle

\abstract{
  Approximating the solutions of boundary value problems governed by partial differential equations with neural networks is challenging, largely due to the difficult training process. This difficulty can be partly explained by the spectral bias, that is, the slower convergence of high-frequency components, and can be mitigated by localizing neural networks via (overlapping) domain decomposition. We combine this localization with the Gauss--Newton method as the optimizer to obtain faster convergence than gradient-based schemes such as Adam; this comes at the cost of solving an ill-conditioned linear system in each iteration. Domain decomposition induces a block-sparse structure in the otherwise dense Gauss--Newton system, reducing the computational cost per iteration. Our numerical results indicate that combining localization and Gauss--Newton optimization is promising for neural network-based solvers for partial differential equations.
}

\section{Introduction}

In recent years, the use of neural networks (NNs) for solving boundary value problems governed by partial differential equations, as an alternative to classical discretizations such as finite differences or finite elements, has been explored. The approach dates back to the 1990s~\cite{dissanayake_neural-network-based_1994,lagaris_artificial_1998}, with popular recent methods including physics-informed NNs (PINNs)~\cite{raissi_physics-informed_2019} and the deep Ritz method~\cite{e_deep_2018}. Implementations thrive on state-of-the-art deep learning frameworks such as TensorFlow, PyTorch, or JAX. Commonly cited additional advantages include their mesh-free, nonlinear nature, improved scaling to higher-dimensional problems, and applicability to parametrized and inverse settings. However, despite these benefits, hyperparameter selection and training remain sensitive, and achieving high accuracy reliably is challenging, especially for complex physical problems.

Difficulties in training NNs are often attributed to the spectral bias (also called the frequency principle), meaning that networks learn low-frequency modes faster than high-frequency modes; cf.~\cite{rahaman_spectral_2019}. This phenomenon can also be linked to the use of activation functions with global support; cf.~\cite{hong_activation_2022}. Consequently, high-frequency errors are reduced slowly, which is particularly problematic for problems with high- and multi-frequency solutions. Among other approaches, localizing the network may help alleviate this issue. For example, finite-basis PINNs (FBPINNs)~\cite{moseley_finite_2023} employ window functions based on an overlapping domain decomposition (DD) to localize networks to subdomains. This strategy improves both convergence and accuracy. The multilevel extension~\cite{dolean_multilevel_2024} further enhances the performance on high- and multi-frequency problems, and extensions to time-dependent problems and operator learning~\cite{heinlein_multifidelity_2025}, randomized neural networks~\cite{anderson_elm-fbpinn_2024,shang_overlapping_2025}, and Kolmogorov--Arnold networks~\cite{howard_finite_2024} have also been explored.

Another contributing factor is the widespread use of first-order optimizers, such as stochastic gradient descent (SGD) or Adam~\cite{kingma_adam_2017}. An alternative is based on the Gauss--Newton (GN) method; cf.~\cite{cai_structure-guided_2024}. For partial differential equation-based loss functions, this is also called the energy natural gradient (ENG) method~\cite{muller_achieving_2023}. This approach can improve convergence, at the cost of solving a linear system involving the Gramian matrix, which may be (nearly) singular and generally ill-conditioned.

In this work, we propose a framework that combines the one-level FBPINN approach with GN-based training, achieving high accuracy with faster convergence. Because of localization in the FBPINN architecture, we obtain a block-sparse GN system, with off-diagonal blocks between non-overlapping subdomains equal to zero; the structure is analogous to that observed for DD-based randomized NNs in~\cite{shang_overlapping_2025}. This can reduce the computational cost of a single GN iteration.

\section{Finite-basis physics-informed neural networks}

For simplicity, we introduce the NN-based approach for a Laplace problem on a Lipschitz domain $\Omega \subset \mathbb{R}^d$ with boundary $\partial \Omega$:
\begin{equation} \label{eq:laplace}
  -\Delta u = f \quad \text{in } \Omega,
  \qquad
  u = g \quad \text{on } \partial\Omega.
\end{equation}
Here, $f$ and $g$ are the given source term and boundary data, respectively. Since we employ PINNs, which rely on the strong form of the PDE, we assume that $u \in \mathcal{C}^2(\Omega)$. Other types of PDEs and boundary conditions can be treated analogously.

\subsection{Physics-informed neural networks}

In order to approximate the solution of~\cref{eq:laplace} using NNs, we introduce the equivalent minimization problem
$$
  \arg\min_{u \in \mathcal{C}^2(\Omega)}
  \| \Delta u + f \|_{L^2(\Omega)}^2
  +
  \omega
  \| u - g \|_{L^2(\partial\Omega)}^2,
$$
with weighting parameter $\omega > 0$. Then, we approximate the integrals in the $L^2$ norms using sampling and approximate the solution using a neural network $u_\theta$, where $\theta \in \mathbb{R}^P$ are the trainable parameters of the NN. We consider simple $k$-layer dense feedforward NNs of the form
$$
  u_\theta(x)
  =
  (L^{(k)} \circ \sigma \circ L^{(k-1)} \circ \sigma \circ \cdots \circ \sigma \circ L^{(1)})(x),
$$
where the  $L^{(i)}(x) = W^{(i)} x + b^{(i)}$ are affine linear transformations with weight matrices $W^{(i)} \in \mathbb{R}^{d_{i} \times d_{i-1}}$ and bias vectors $b^{(i)} \in \mathbb{R}^{d_{i}}$, and $\sigma$ is a nonlinear activation function applied component-wise. In order to enforce boundary conditions, we introduce constraint operator $\mathcal{C}$, which enforces them  explicitly; as a result, the boundary condition loss vanishes. This results in the discrete minimization problem

\begin{figure}[t]
  \centering
  \scalebox{0.8}{
  \begin{tikzpicture}
    \pgfmathsetmacro{\k}{5}         
    \pgfmathsetmacro{\delta}{0.5}   
    \pgfmathsetmacro{\ramp}{2*\delta/\k}  
    \pgfmathsetmacro{\xjsZero}{-1 - 2*\delta}
    \pgfmathsetmacro{\xjeZero}{(2*(0+1)+\delta)/\k - 1}
    \pgfmathsetmacro{\xjsOne}{(2*(1)-\delta)/\k - 1}
    \pgfmathsetmacro{\xjeOne}{(2*(1+1)+\delta)/\k - 1}
    \pgfmathsetmacro{\xjsTwo}{(2*(2)-\delta)/\k - 1}
    \pgfmathsetmacro{\xjeTwo}{(2*(2+1)+\delta)/\k - 1}
    \pgfmathsetmacro{\xjsThree}{(2*(3)-\delta)/\k - 1}
    \pgfmathsetmacro{\xjeThree}{(2*(3+1)+\delta)/\k - 1}
    \pgfmathsetmacro{\xjsFour}{(2*(4)-\delta)/\k - 1}
    \pgfmathsetmacro{\xjeFour}{1 + 2*\delta}

    \begin{axis}[
        width=\textwidth,
        height=0.48\textwidth,
        domain=-1.0:1.0,
        samples=400,
        xmin=-1, xmax=1,
        ymin=0, ymax=1.1,
        grid=both,
        xlabel={$x$}, ylabel={$\omega_j(x)$},
        legend style={at={(0.5,1.02)}, anchor=south, legend columns=5},
        tick label style={/pgf/number format/fixed},
      ]


      \addplot [blue, thick]
      { ( (x<\xjsZero) || (x>=\xjeZero) ) ? 0 :
        ( (x<\xjsZero+\ramp) ?
        0.5*(1 - cos(deg(pi*(\k)*(x-\xjsZero)/(2*\delta)))) :
        ( (x<\xjeZero-\ramp) ? 1 :
        ( (x<\xjeZero) ?
        0.5*(1 - cos(deg(pi*(\k)*(x-\xjeZero)/(2*\delta)))) : 0 ) ) ) };
      \addlegendentry{$\omega_0$}

      \addplot [red, thick]
      { ( (x<\xjsOne) || (x>=\xjeOne) ) ? 0 :
        ( (x<\xjsOne+\ramp) ?
        0.5*(1 - cos(deg(pi*(\k)*(x-\xjsOne)/(2*\delta)))) :
        ( (x<\xjeOne-\ramp) ? 1 :
        ( (x<\xjeOne) ?
        0.5*(1 - cos(deg(pi*(\k)*(x-\xjeOne)/(2*\delta)))) : 0 ) ) ) };
      \addlegendentry{$\omega_1$}

      \addplot [teal!70!black, thick]
      { ( (x<\xjsTwo) || (x>=\xjeTwo) ) ? 0 :
        ( (x<\xjsTwo+\ramp) ?
        0.5*(1 - cos(deg(pi*(\k)*(x-\xjsTwo)/(2*\delta)))) :
        ( (x<\xjeTwo-\ramp) ? 1 :
        ( (x<\xjeTwo) ?
        0.5*(1 - cos(deg(pi*(\k)*(x-\xjeTwo)/(2*\delta)))) : 0 ) ) ) };
      \addlegendentry{$\omega_2$}

      \addplot [orange!85!black, thick]
      { ( (x<\xjsThree) || (x>=\xjeThree) ) ? 0 :
        ( (x<\xjsThree+\ramp) ?
        0.5*(1 - cos(deg(pi*(\k)*(x-\xjsThree)/(2*\delta)))) :
        ( (x<\xjeThree-\ramp) ? 1 :
        ( (x<\xjeThree) ?
        0.5*(1 - cos(deg(pi*(\k)*(x-\xjeThree)/(2*\delta)))) : 0 ) ) ) };
      \addlegendentry{$\omega_3$}

      \addplot [purple, thick]
      { ( (x<\xjsFour) || (x>=\xjeFour) ) ? 0 :
        ( (x<\xjsFour+\ramp) ?
        0.5*(1 - cos(deg(pi*(\k)*(x-\xjsFour)/(2*\delta)))) :
        ( (x<\xjeFour-\ramp) ? 1 :
        ( (x<\xjeFour) ?
        0.5*(1 - cos(deg(pi*(\k)*(x-\xjeFour)/(2*\delta)))) : 0 ) ) ) };
      \addlegendentry{$\omega_4$}

      \addplot [black, dashed] coordinates {(-1,1) (1,1)};
      \addlegendentry{$\sum_j \omega_j(x) \approx 1$}
    \end{axis}
  \end{tikzpicture}
  }
  \caption{Five overlapping partition of unity window functions $\omega_k$ on $[-1, 1]$ built from cosine (Hann) windows on overlapping subintervals.}
  \label{fig:pou}
\end{figure}
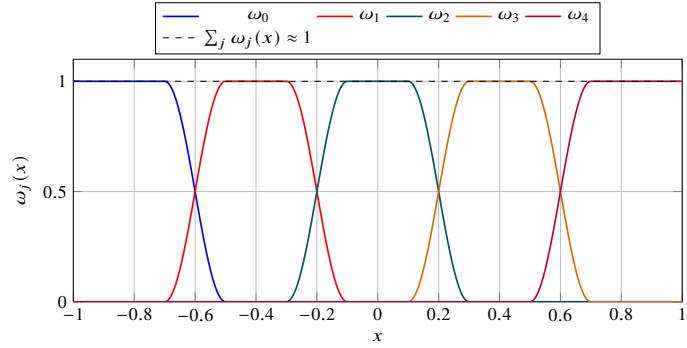

\begin{equation} \label{eq:loss}
  \arg\min_{\theta} \underbrace{\frac{1}{N} \sum_{i=1}^N \left( \Delta \mathcal{C} (u_\theta(x_i)) + f(x_i) \right)^2}_{\eqqcolon \mathcal{L}}
\end{equation}
The points $\{x_i\}_{i=1}^N \subset \Omega$ are sampled collocation points in the domain. There are various choices for the sampling; see~\cite{wu_comprehensive_2023} for more details. To solve~\cref{eq:loss}, standard gradient-based optimizers, such as Adam~\cite{kingma_adam_2017}, are often employed.

As mentioned before, training is particularly challenging when the solution exhibits high-frequency components, due to the spectral bias~\cite{rahaman_spectral_2019}; this is further amplified by the spectral properties of the differential operator.

\subsection{Domain decomposition-based network architecture}

\begin{table}[t]
  \centering
  \caption{\textit{Test case 1:} hyperparameters for the domain-decomposed PINN.}
  \begin{tabular}{l@{\hspace{4mm}}l@{\hspace{4mm}}l}
    \toprule
    \textbf{optimizer}  & \textbf{Adam}                                                                               & \textbf{Gauss--Newton}            \\
    learning rate       & $lr = 10^{-2}$                                                                              & constant step size $\eta=10^{-2}$ \\
    stopping criterion  & $2000$ its.\ or $\mathcal{L} < 10^{-6}$                                                             & $5000$ its.\ or $\mathcal{L} < 10^{-6}$   \\
    \midrule
    collocation points  & \multicolumn{2}{l}{$N_f=1000$ (uniform)}                                                                                        \\
    subdomains          & \multicolumn{2}{l}{ $k=24$; overlap $\delta=0.5$ }                                                                              \\
    constraint operator & \multicolumn{2}{l}{$ \mathcal{C}(x) = \tanh(k\pi x)$}                                                                              \\
    subdomain network   & \multicolumn{2}{l}{MLP $[1,20,1]$, $\sigma = \tanh$; init.: $W\sim\mathcal U[-1,1]$, $b=0$}                                     \\
    \bottomrule
  \end{tabular}
  \label{Tab1}
\end{table}

\begin{table}[t]
  \centering
  \caption{\textit{Test case 2:} hyperparameters for the 2D domain-decomposed PINN.}
  \begin{tabular}{l@{\hspace{4mm}}l@{\hspace{4mm}}l}
    \toprule
    \textbf{optimizer}  & \textbf{Adam}                                                                  & \textbf{Gauss--Newton}                              \\
    learning rate       & $10^{-3}$                                                                      & constant step size $\eta=10^{-2}$                   \\
    stopping criterion  & $\num{30000}$ its.\ or $\mathcal L<10^{-5}$                           & $\num{1000}$ its.\ or $\mathcal L<10^{-5}$ \\
    \midrule
    collocation points  & \multicolumn{2}{l}{$100{\times}100$ (uniform grid) $\Rightarrow N_f=\num{10000}$}                                                       \\
    subdomains          & \multicolumn{2}{l}{$(k_x,k_y)=(2,2)$; overlap $(\delta_x,\delta_y)=(0.5,0.5)$}                                                       \\
    constraint operator & \multicolumn{2}{l}{$\mathcal{C}(x) = \phi(x)\phi(y)$, where $\phi(s)=1-s^2$}                                                                    \\
    subdomain network   & \multicolumn{2}{l}{MLP $[2,20,1]$, $\sigma = \tanh$; init.: $W\sim\mathcal U[-1,1]$, $b=0$}                                                        \\
    \bottomrule
  \end{tabular}
  \label{Tab2}
\end{table}

In order to localize the neural networks and improve the approximation of high-frequency components, we employ the one-level FBPINN approach~\cite{moseley_finite_2023}.
Specifically, we introduce a set of $K$ overlapping subdomains
$\{\Omega_k\}_{k=1}^K$ such that
\[
  \bigcup_{k=1}^K \Omega_k = \Omega,
\]
and a corresponding set of window functions $\{\omega_k\}_{k=1}^K$ forming a partition of unity, that is,
\[
  \sum_{k=1}^K \omega_k(x) = 1
  \quad \text{for all } x \in \Omega,
  \qquad
  \omega_k(x) = 0
  \quad \text{for } x \in \Omega \setminus \Omega_k.
\]

We define the overlap between neighboring subdomains by a the parameter $\delta > 0$. For simplicity, let us consider the one-dimensional case; square and cubic domains on two and three dimensions follow via tensor product. With
$r = \frac{2\delta}{K}$, we define the window function $\omega_k(x)$ for the $k$-th subdomain $\Omega_k = (a_k, b_k)$ as
\[
\omega_k(x) =
\begin{cases}
1/2 \!\left[1 - \cos\!\left(\pi (x - a_k) / r\right)\right],
  & a_k \le x < a_k + r, \\
1, & a_k + r \le x < b_k - r, \\
1/2 \!\left[1 - \cos\!\left(\pi (x - b_k)/r\right)\right],
  & b_k - r \le x < b_k, \\
0, & \text{otherwise.}
\end{cases}
\]
with
\[
a_k = (2k - \delta) / K - 1,
\qquad
b_k = (2(k+1) + \delta) /K - 1,
\qquad
k = 0, 1, \dots, K-1.
\]

For each overlapping subdomain, we introduce a neural network $u_{\theta_k}$ and construct the global solution as
\begin{equation} \label{eq:fbpinn}
  u_\theta(x) = \sum_{k=1}^K \omega_k(x)\, u_{\theta_k}\bigl(n_k(x)\bigr),
\end{equation}
where $\theta = (\theta_1, \ldots, \theta_K)$ are the trainable parameters of the local networks, and $n_k(x)$ is a normalization function mapping $\Omega_k$ to a reference domain; see~\cite{moseley_finite_2023} for details. Finally, we replace $u_\theta$ in~\cref{eq:loss} with the FBPINN ansatz~\cref{eq:fbpinn}.

The use of window functions and normalization allows each subdomain network to capture locally high-frequency components efficiently.
For problems exhibiting multiple frequency scales, a multilevel extension of the FBPINN framework~\cite{dolean_multilevel_2024} can further enhance performance.

\section{Gauss--Newton training}

Usually, the minimization problem in~\cref{eq:loss} is optimized using methods based on gradient descent,
$$
  \theta^{(m+1)}
  =
  \theta^{(m)} - \alpha \nabla_\theta \mathcal{L}(\theta^{(m)}),
$$
where $\theta^{(m)}$ are the parameters at iteration $m$ and $\alpha$ is the learning rate. The most popular extension of gradient descent is the Adam optimizer~\cite{kingma_adam_2017}, which employs adaptive learning rates and momentum terms.

As discussed in~\cite{cai_structure-guided_2024,muller_achieving_2023}, the training performance can be significantly improved by employing the GN method, which yields the following update:
\begin{equation} \label{eq:gn-update}
  \theta^{(m+1)}
  =
  \theta^{(m)} - \alpha\, G ^+ (\theta^{(m)})\, \nabla_{\theta} \mathcal{L}(\theta^{(m)}),
\end{equation}
where $G^+(\theta)$ is the pseudoinverse of the energy Gram matrix. For the PINN loss~\cref{eq:loss}, it reads
$$
  G(\theta)_{ij}
  =
  \int_\Omega \Delta (\partial_{\theta_i} (\mathcal{C} u_\theta)) \, \Delta (\partial_{\theta_j} (\mathcal{C} u_\theta)) \, dx.
$$
GN can be understood as approximating the Hessian of Newton's method by the energy Gram matrix. Note that $G(\theta)$ can equivalently be expressed as $J(\theta)^T J(\theta)$, where $J(\theta)$ denotes the Jacobian of the PINN loss. The term
$G^+(\theta)\, \nabla_{\theta} \mathcal{L}(\theta)$,
appearing in the update~\cref{eq:gn-update}, is also called the energy natural gradient (ENG)~\cite{muller_achieving_2023}.

Because $G(\theta)$ can be nearly singular, we regularize it as $G(\theta) + \mu I$ with a small $\mu > 0$. For this work, $\mu=1$ is chosen for all numerical experiments.

\section{Numerical results}

\begin{table}[t]
  \centering
  \caption{Relative $\ell_2$ test errors of FBPINN with Adam and Gauss--Newton training for the 1D and 2D problems.
  }
  \begin{tabular}{l@{\hspace{4mm}}l@{\hspace{4mm}}l}
    \toprule
    \textbf{Problem} & \textbf{FBPINN (Adam)} & \textbf{FBPINN (Gauss--Newton)} \\
    \midrule
    Test case 1: 1D ODE           & 7.8 $\times$ $10^{-3}$                & 8.0 $\times$ $10^{-4}$        \\
    Test case 2: 2D Helmholtz     & 2.0 $\times$ $10^{-4}$                & 1.5 $\times$ $10^{-4}$                         \\
    \bottomrule
  \end{tabular}
  \label{tab:rel_errors}
\end{table}

\begin{figure}[t]
  \centering

  \begin{minipage}{0.31\textwidth}
    \centering

    \hspace*{3mm} PINN \\[3pt]
    \includegraphics[height=0.7\textwidth]{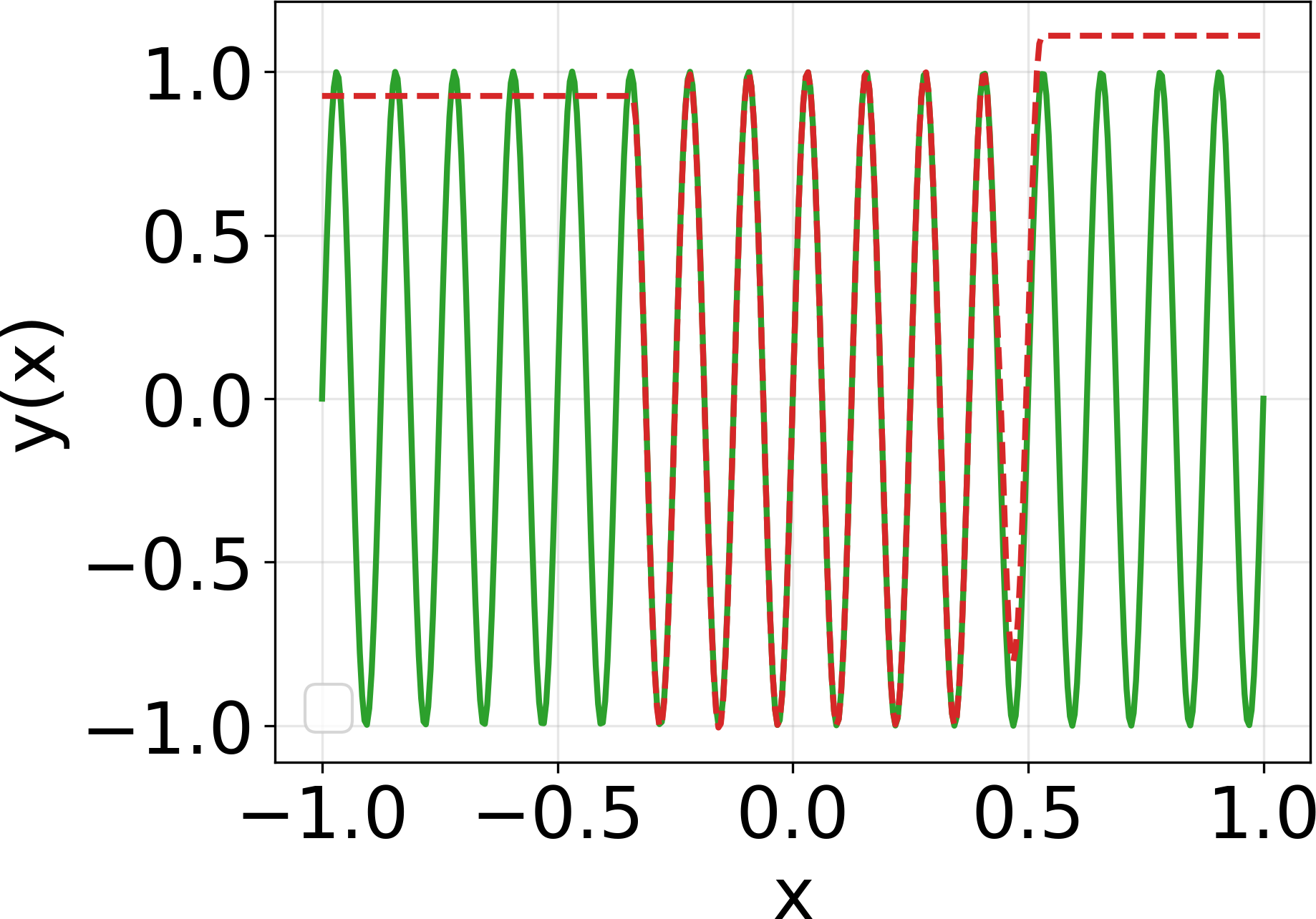} \\
    \includegraphics[height=0.7\textwidth]{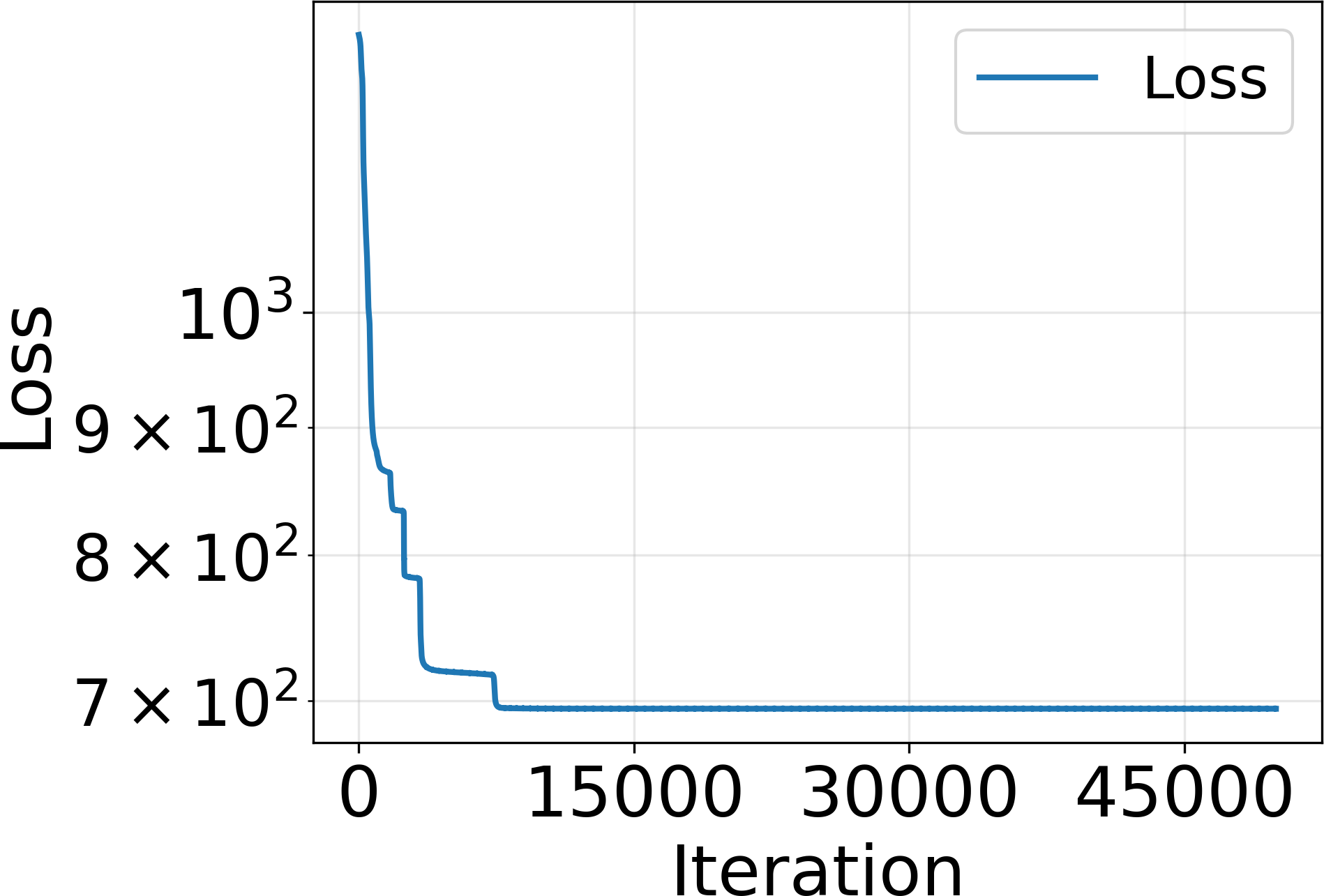}

  \end{minipage}
  \hspace*{0.01\textwidth}
  \begin{minipage}{0.31\textwidth}
    \centering

    \hspace*{2mm} FBPINN (Adam) \\[3pt]
    \includegraphics[height=0.7\textwidth]{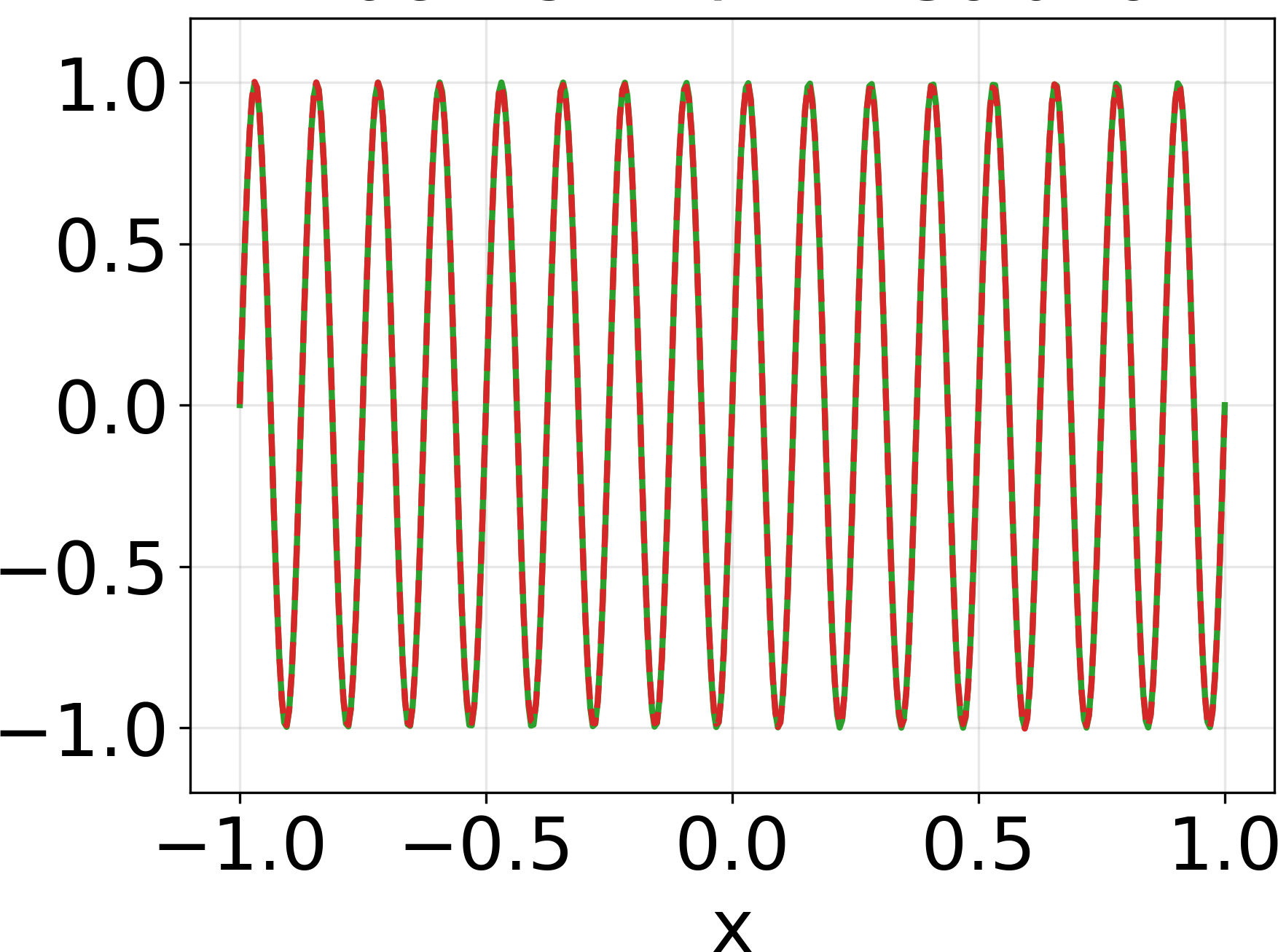} \\
    \includegraphics[height=0.7\textwidth]{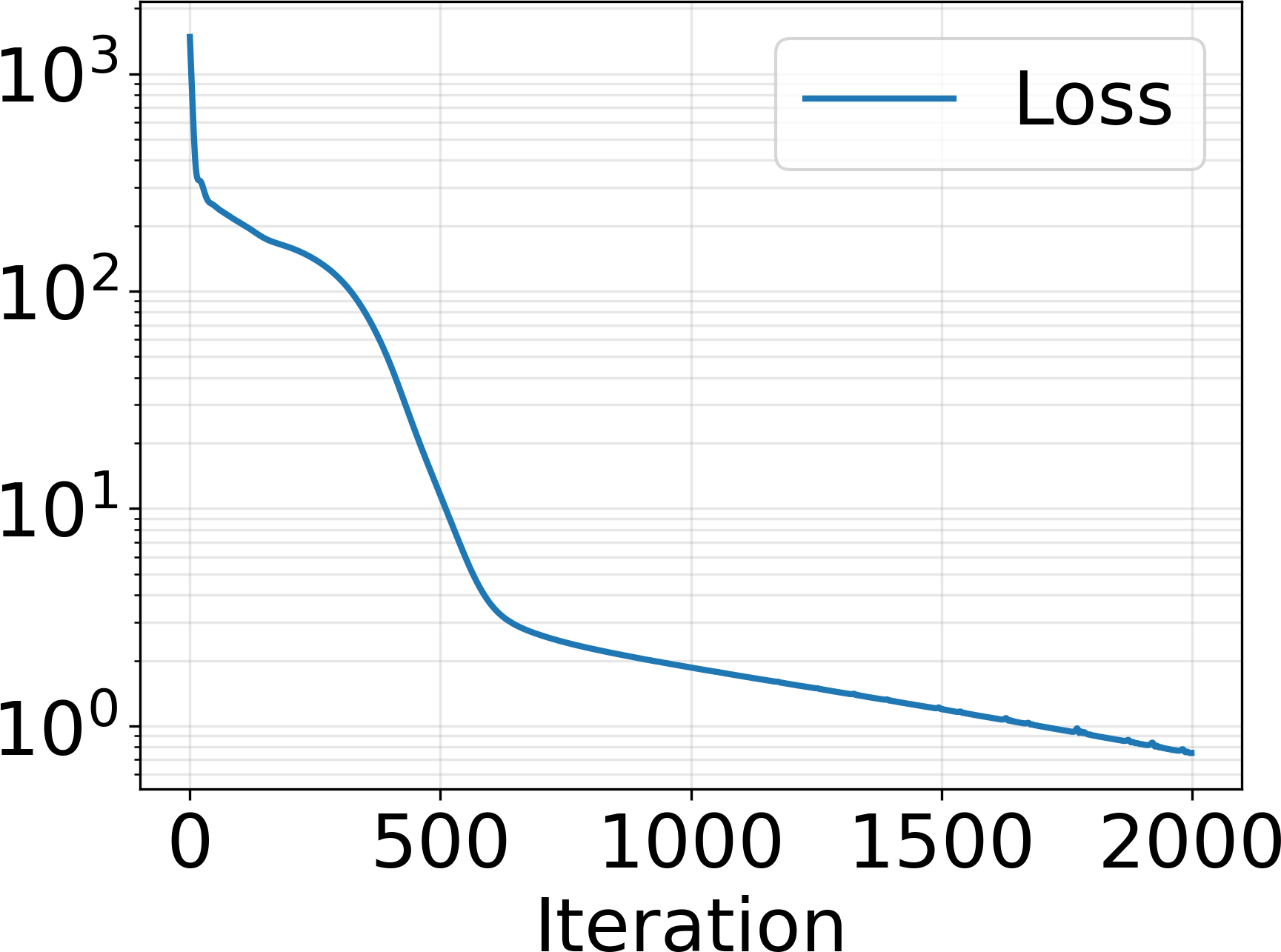}

  \end{minipage}
  \begin{minipage}{0.31\textwidth}
    \centering

    FBPINN (Gauss--Newton) \\[3pt]
    \includegraphics[height=0.7\textwidth]{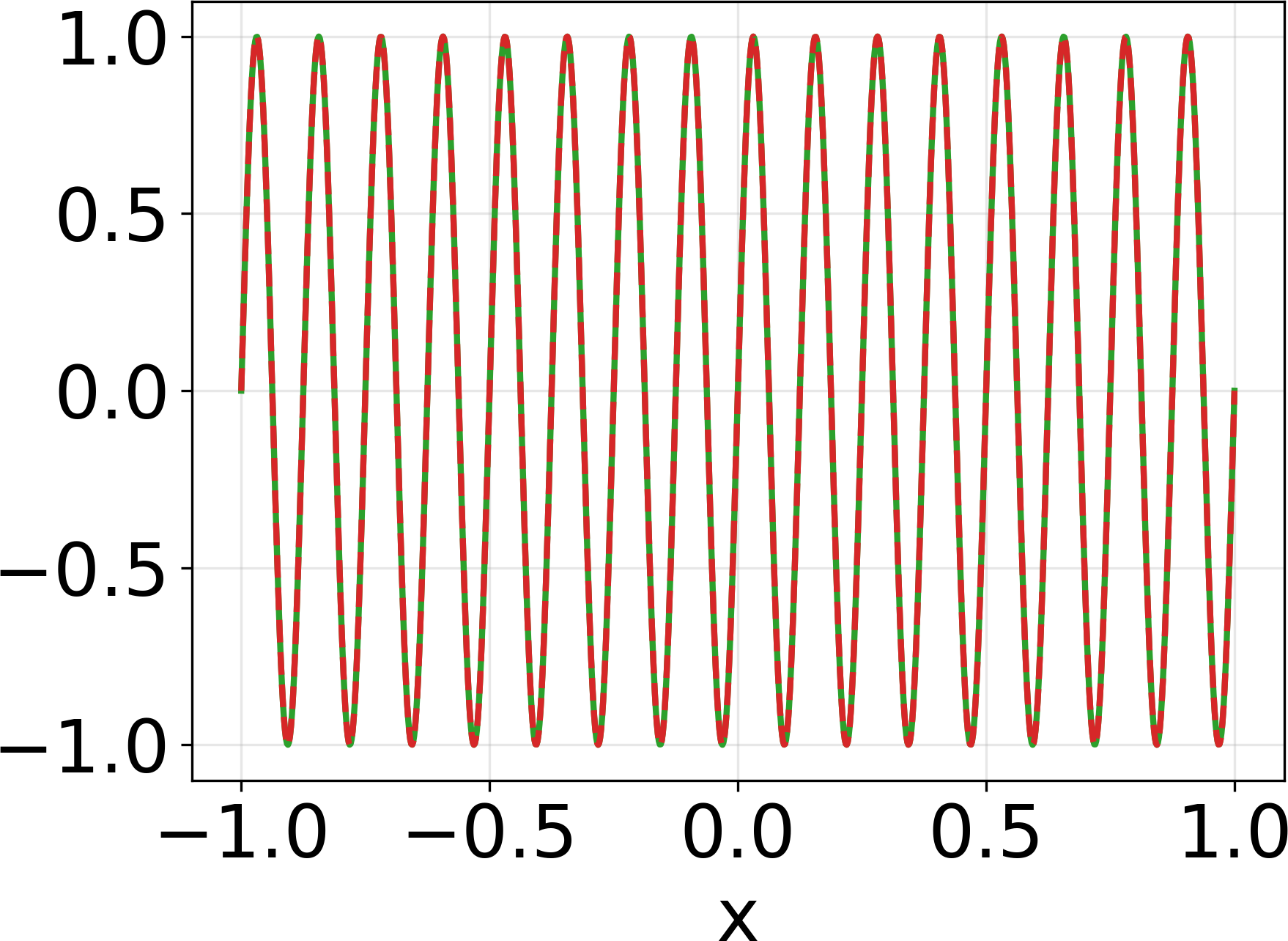} \\
  \includegraphics[height=0.7\textwidth]{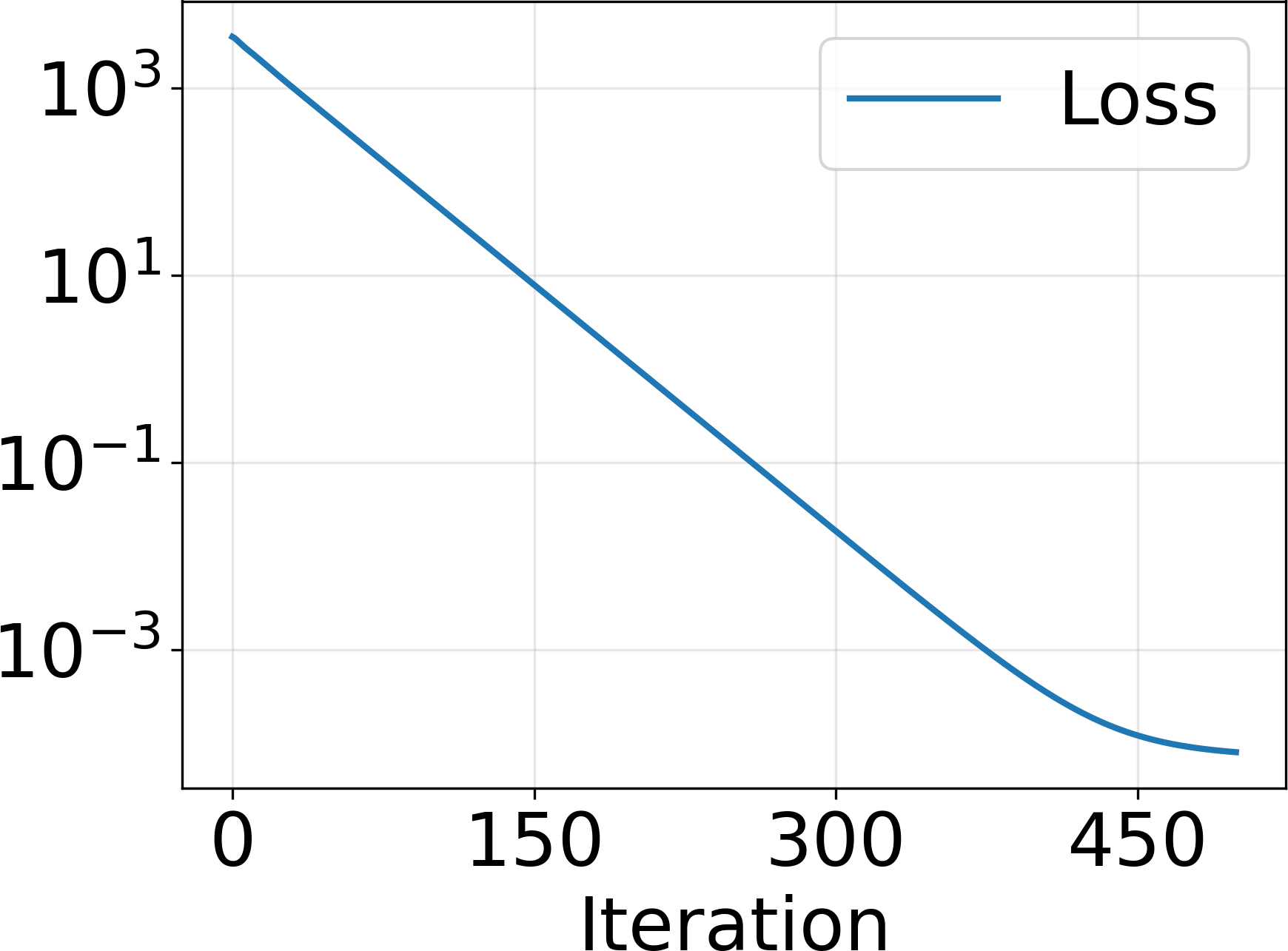}

  \end{minipage}
  \caption{
    Comparison of prediction in orange against the analytical solution in green (top) and training loss (bottom) of three methods for the ordinary differential equation problem: (left) vanilla PINN, (middle) FBPINN, and (right) Gauss--Newton PINN. Left panels show the loss on a logarithmic scale versus iteration; right panels overlay the learned solution with the exact solution. The green curve shows the exact solution, and the red curve shows the prediction.
  }
  \label{fig2}
\end{figure}

In this section, we validate the performance of the proposed approach using two canonical problems. The first test case is
\[
  \frac{du}{dx} - 16\pi \cos\bigl(16\pi x\bigr) = 0, \qquad x \in [-1, 1],
\]
with $u(0) = 0$. As the second case, we consider the two-dimensional Helmholtz equation
\[
  \Delta u(x,y) + k^{2} u(x,y) = f(x,y), \quad (x,y) \in [-1,1]^2,
\]
subject to homogeneous Dirichlet boundary conditions. We consider a low wave number case with \(k = 1\) and the source term
\[
  f(x,y) = -4 + 2\bigl(x^{2} + y^{2}\bigr) + \bigl(1 - x^{2}\bigr)\bigl(1 - y^{2}\bigr),
\]
which yields the exact solution
$
  u_{\text{exact}}(x,y) = \bigl(1 - x^{2}\bigr)\bigl(1 - y^{2}\bigr).
$

As the baseline for the 1D problem, we employ a fully connected NN with three hidden layers of 20 neurons each, a $\tanh$ activation, and Glorot initialization~\cite{glorot2010understanding}. We place $N_f = 256$ collocation points uniformly in $[-1,1]$ and train with the Adam optimizer~\cite{kingma_adam_2017} using a learning rate of $10^{-3}$ and full batches for \num{50000} steps. This vanilla PINN attains a relative error of $1.31$ on test case~1. We compare this baseline with FBPINNs trained using Adam and Gauss--Newton, with hyperparameters summarized in~\Cref{Tab1}.

\begin{figure}[t]
  \centering

  \includegraphics[height=0.34\textwidth]{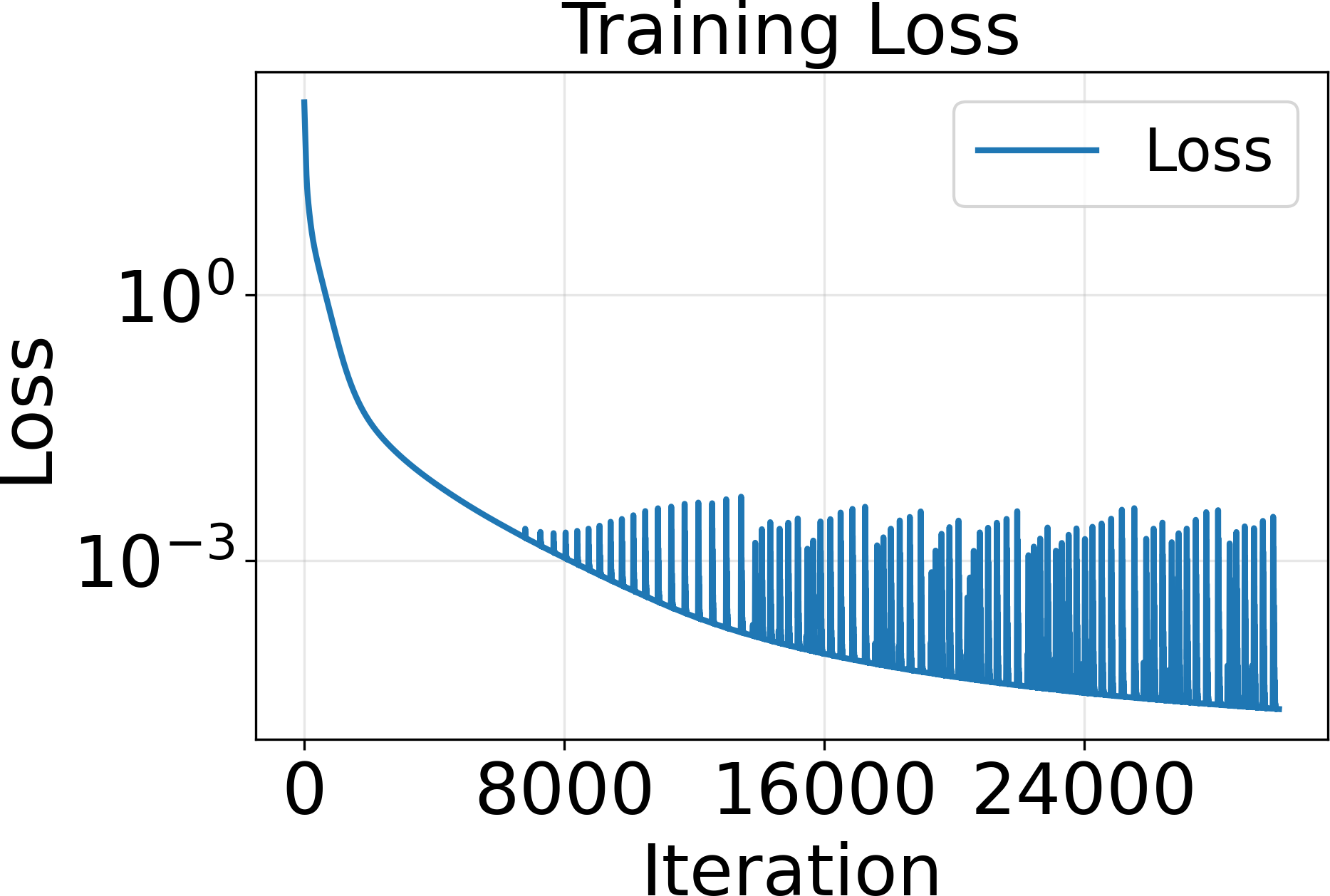}
  \hspace*{0.01\textwidth}
  \includegraphics[height=0.34\textwidth]{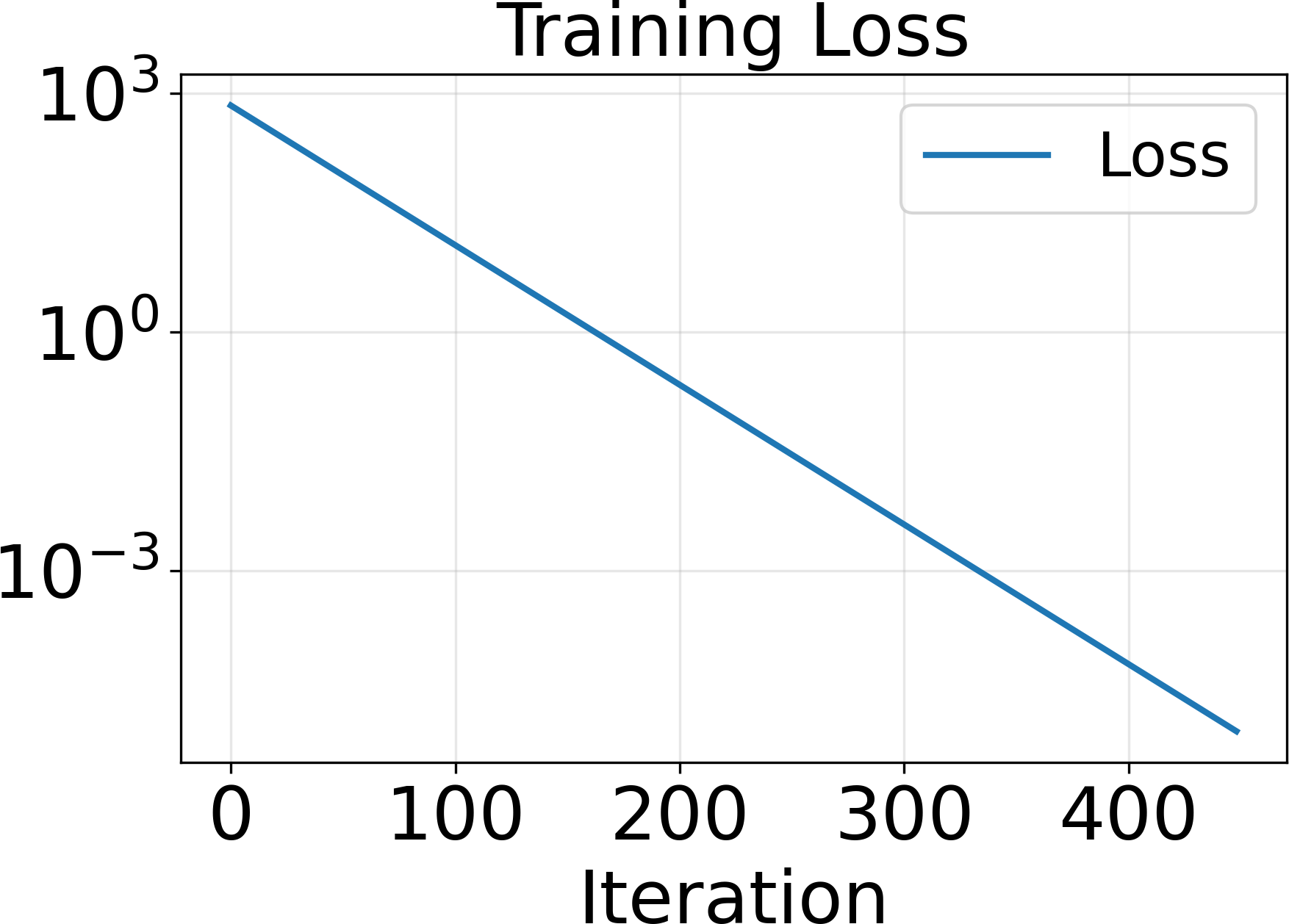}
  \caption{Training loss for the 2D Helmholtz problem: Adam (left) and Gauss--Newton (right).}
  \label{fig3}
\end{figure}

\begin{figure}[t]
  \centering

  \includegraphics[height=0.22\textwidth]{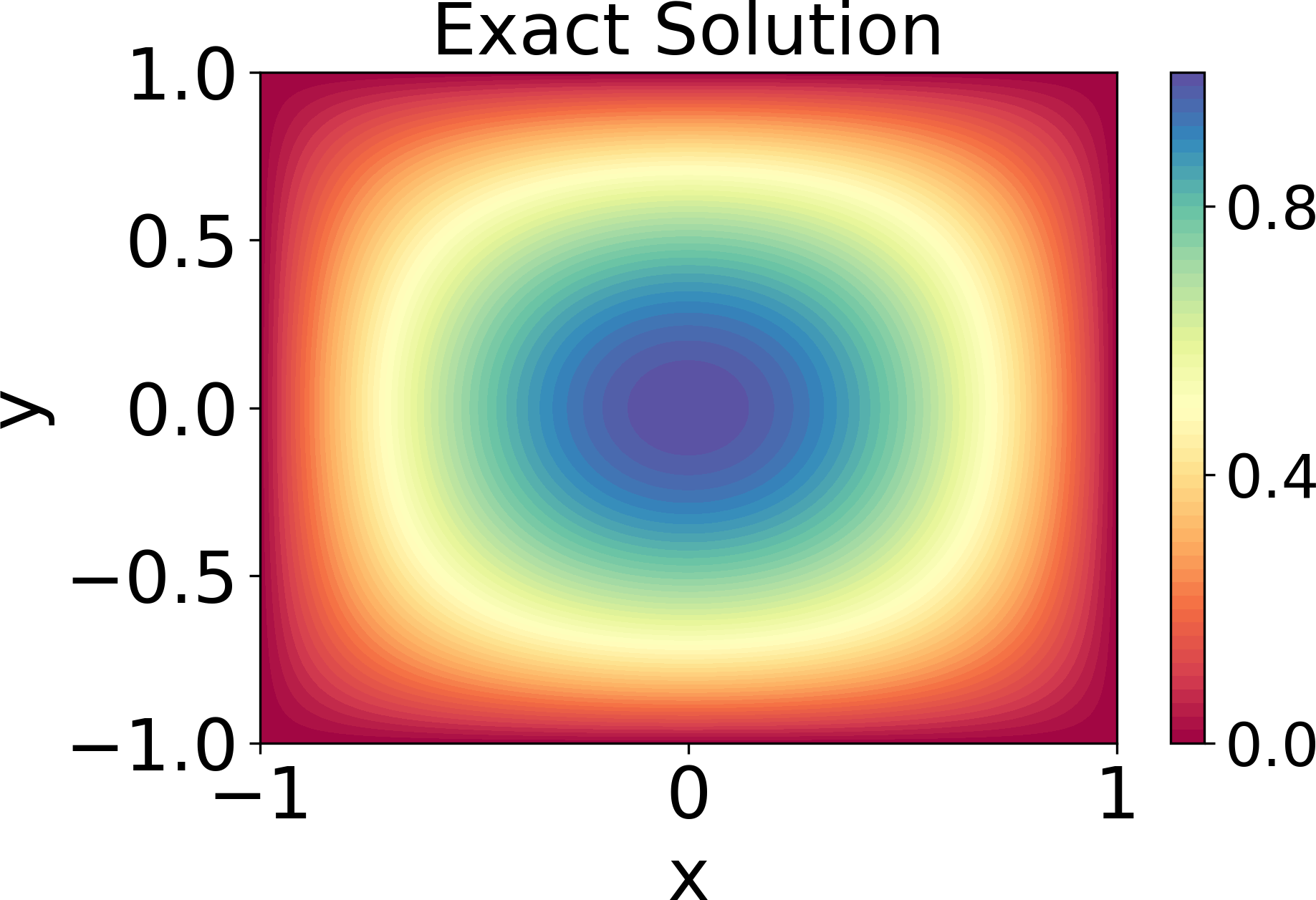}
  \hspace*{0.01\textwidth}
  \includegraphics[height=0.22\textwidth]{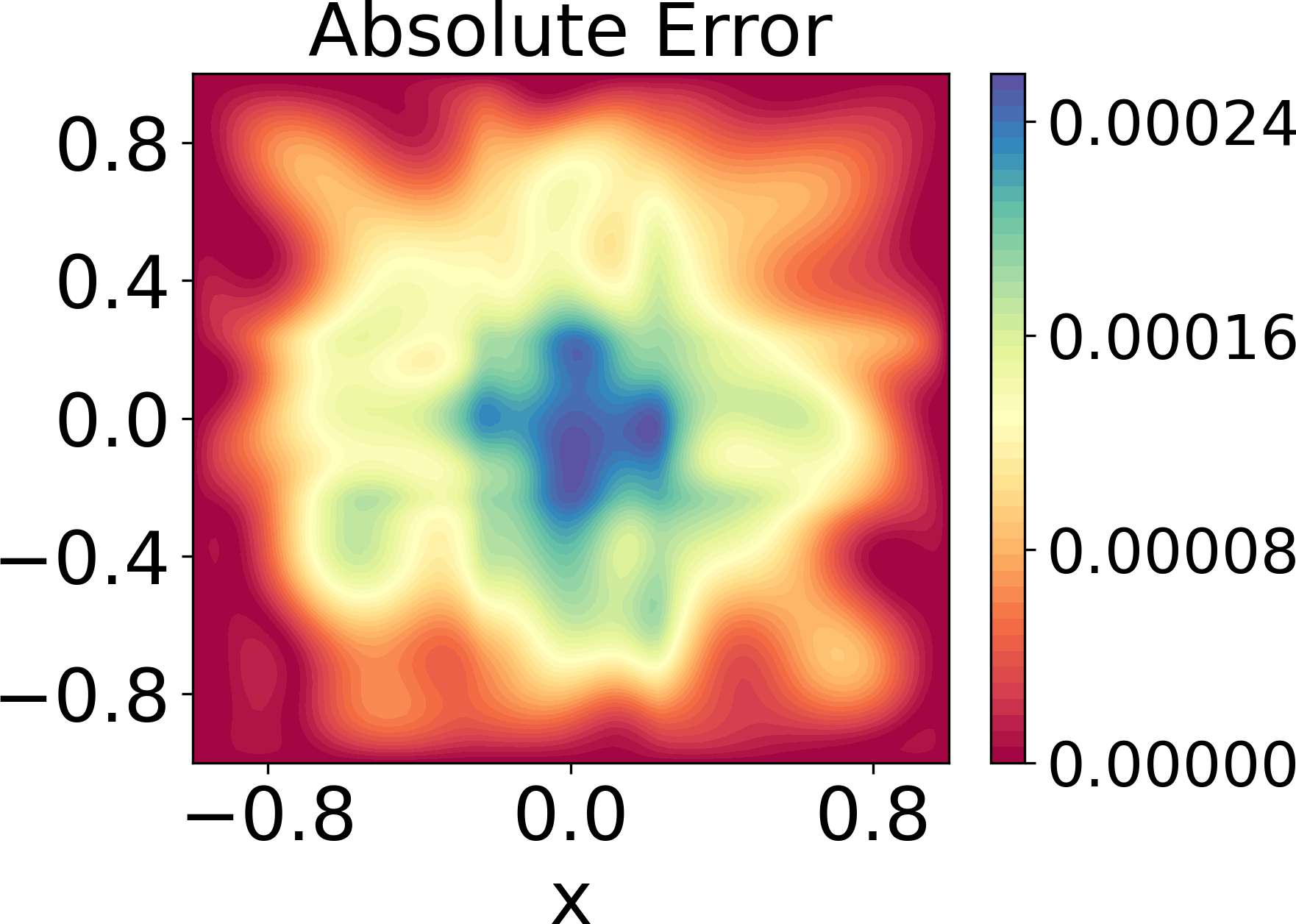}
  \hspace*{0.01\textwidth}
  \includegraphics[height=0.22\textwidth]{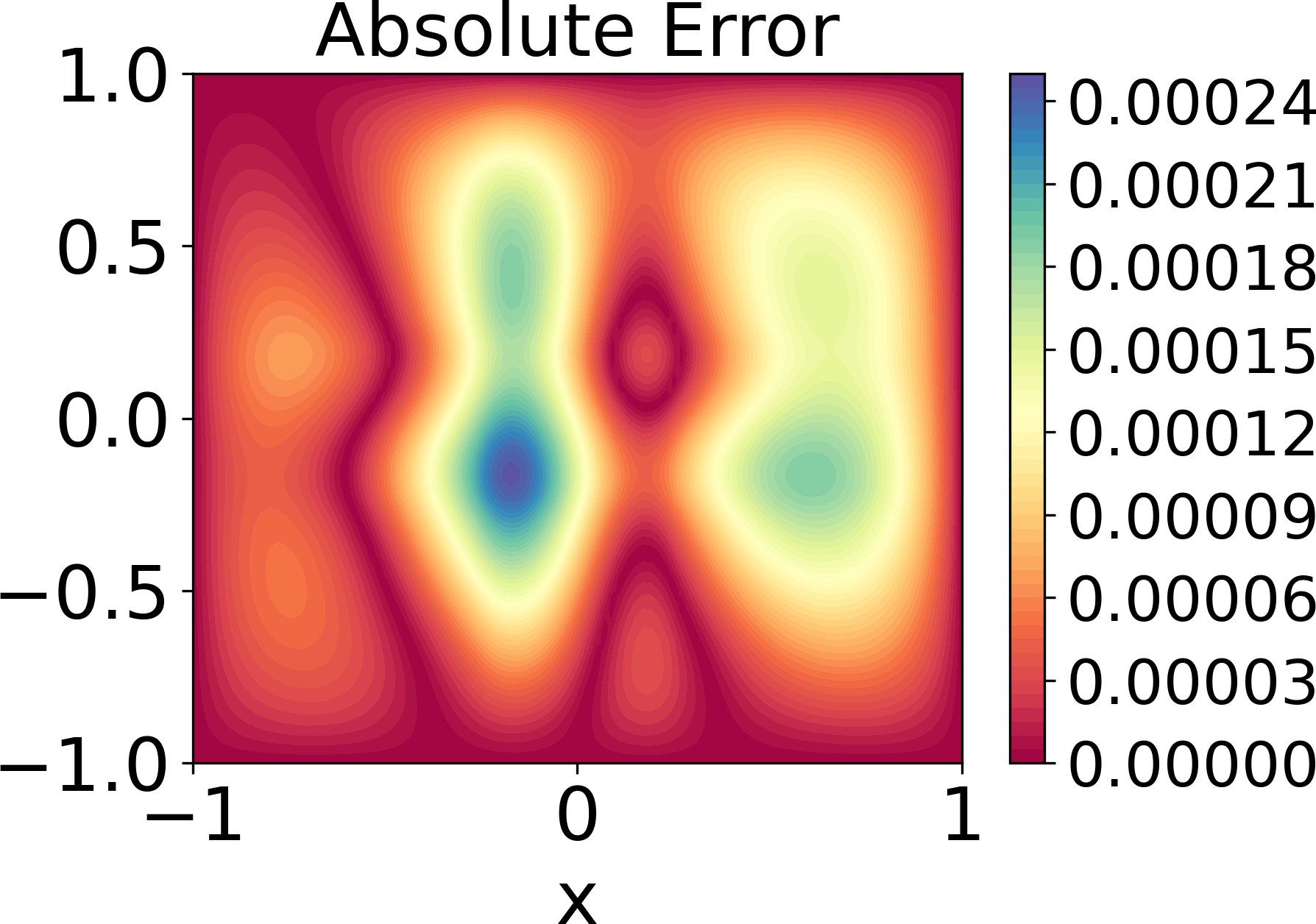}
  \caption{Exact solution (left) and absolute errors for Adam (middle) and Gauss--Newton (right).}
  \label{fig4}
\end{figure}

\Cref{tab:rel_errors} reports the relative $\ell_2$ errors for the 1D ODE. FBPINN trained with Adam attains a test error of $7.8\times 10^{-3}$, improving over the baseline PINN, and Gauss--Newton reduces the error by an order of magnitude to $8.0\times 10^{-4}$. The loss and prediction curves in~\Cref{fig2} show that the vanilla PINN converges slowly, the FBPINN improves accuracy, and Gauss--Newton further accelerates convergence while aligning the prediction with the reference solution.

\begin{wrapfigure}[14]{r}{0.3\textwidth}
  \centering

  \vspace*{-6.5mm}

  \includegraphics[width=0.3\textwidth]{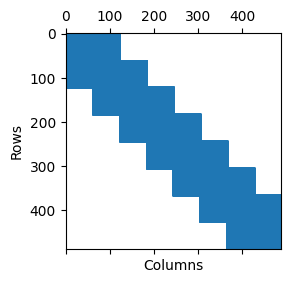}
  \caption{Exemplary sparisity pattern of the Gramian $G$ for the FBPINN architecture for the one dimensional ODE problem with $8$ subdomains.}
  \label{fig:sparsity}
\end{wrapfigure}

We did not conduct a detailed study on the NN architectures. Nonetheless, even PINN models with more parameters continue to suffer from spectral bias and difficulties in learning highly oscillatory functions; see, for example,~\cite{moseley_finite_2023,dolean_multilevel_2024,kapoor2024transfer}.

While we do not report computing times, \Cref{fig:sparsity} illustrates the sparsity structure of the Gramian $G$ induced by the FBPINN architecture. A fully connected network would yield a dense Gramian, but the domain decomposition creates pronounced sparsity, allowing for using efficient iterative solvers such as those in~\cite{shang_overlapping_2025}.

The 2D Helmholtz problem shows the same pattern. We compare FBPINNs trained with Adam and Gauss--Newton using the hyperparameters in~\Cref{Tab2}. The relative errors fall from $2.0\times 10^{-4}$ with Adam to $1.5\times 10^{-4}$ with Gauss--Newton, confirming a consistent improvement; see the error plots in~\Cref{fig4}. The loss curves in~\Cref{fig3} likewise show slower convergence for Adam and faster decay with closer agreement for Gauss--Newton.

Overall, training with the Gauss--Newton method significantly accelerates convergence and also improves accuracy across both model problems.

\medskip

\noindent \textbf{Acknowledgments} \ We would like to thank Marius Zeinhofer and Rami Masri for the helpful discussions about the presented approach.

\bibliographystyle{spmpsci}
\bibliography{bibliography}

\begin{thebibliography}{10}
\providecommand{\url}[1]{{#1}}
\providecommand{\urlprefix}{URL }
\expandafter\ifx\csname urlstyle\endcsname\relax
  \providecommand{\doi}[1]{DOI~\discretionary{}{}{}#1}\else
  \providecommand{\doi}{DOI~\discretionary{}{}{}\begingroup \urlstyle{rm}\Url}\fi

\bibitem{anderson_elm-fbpinn_2024}
Anderson, S., Dolean, V., Moseley, B., Pestana, J.: {ELM}-{FBPINN}: efficient finite-basis physics-informed neural networks (2024).
\newblock ArXiv:2409.01949

\bibitem{cai_structure-guided_2024}
Cai, Z., Ding, T., Liu, M., Liu, X., Xia, J.: A {Structure}-{Guided} {Gauss}-{Newton} {Method} for {Shallow} {ReLU} {Neural} {Network} (2024).
\newblock ArXiv:2404.05064

\bibitem{dissanayake_neural-network-based_1994}
Dissanayake, M.W.M.G., Phan-Thien, N.: Neural-network-based approximations for solving partial differential equations.
\newblock Commun. Numer. Methods Eng. \textbf{10}(3), 195--201 (1994)

\bibitem{dolean_multilevel_2024}
Dolean, V., Heinlein, A., Mishra, S., Moseley, B.: Multilevel domain decomposition-based architectures for physics-informed neural networks.
\newblock Comput. Methods Appl. Mech. Eng. \textbf{429}, 117116 (2024)

\bibitem{e_deep_2018}
E, W., Yu, B.: The {Deep} {Ritz} {Method}: {A} {Deep} {Learning}-{Based} {Numerical} {Algorithm} for {Solving} {Variational} {Problems}.
\newblock Commun. Math. Stat. \textbf{6}(1), 1--12 (2018)

\bibitem{glorot2010understanding}
Glorot, X., Bengio, Y.: Understanding the difficulty of training deep feedforward neural networks.
\newblock In: Proceedings of the thirteenth international conference on artificial intelligence and statistics, pp. 249--256. JMLR Workshop and Conference Proceedings (2010)

\bibitem{heinlein_multifidelity_2025}
Heinlein, A., Howard, A.A., Stinis, P., Beecroft, D.: Multifidelity domain decomposition-based physics-informed neural networks and operators for time-dependent problems.
\newblock In: Mathematical {Optimization} for {Machine} {Learning}, pp. 79--92. De Gruyter (2025)

\bibitem{hong_activation_2022}
Hong, Q., Siegel, J.W., Tan, Q., Xu, J.: On the {Activation} {Function} {Dependence} of the {Spectral} {Bias} of {Neural} {Networks} (2022).
\newblock ArXiv:2208.04924

\bibitem{howard_finite_2024}
Howard, A.A., Jacob, B., Murphy, S.H., Heinlein, A., Stinis, P.: Finite basis {Kolmogorov}-{Arnold} networks: domain decomposition for data-driven and physics-informed problems (2024).
\newblock ArXiv:2406.19662

\bibitem{kapoor2024transfer}
Kapoor, T., Wang, H., N{\'u}{\~n}ez, A., Dollevoet, R.: Transfer learning for improved generalizability in causal physics-informed neural networks for beam simulations.
\newblock Eng. Appl. Artif. Intell. \textbf{133}, 108085 (2024)

\bibitem{kingma_adam_2017}
Kingma, D.P., Ba, J.: Adam: {A} {Method} for {Stochastic} {Optimization} (2017).
\newblock ArXiv:1412.6980

\bibitem{lagaris_artificial_1998}
Lagaris, I., Likas, A., Fotiadis, D.: Artificial neural networks for solving ordinary and partial differential equations.
\newblock IEEE Trans. Neural Networks \textbf{9}(5), 987--1000 (1998)

\bibitem{moseley_finite_2023}
Moseley, B., Markham, A., Nissen-Meyer, T.: Finite basis physics-informed neural networks ({FBPINNs}): a scalable domain decomposition approach for solving differential equations.
\newblock Adv. Comput. Math. \textbf{49}(4), 62 (2023)

\bibitem{muller_achieving_2023}
Müller, J., Zeinhofer, M.: Achieving {High} {Accuracy} with {PINNs} via {Energy} {Natural} {Gradient} {Descent}.
\newblock In: Proceedings of the 40th {International} {Conference} on {Machine} {Learning}, pp. 25471--25485. PMLR (2023)

\bibitem{rahaman_spectral_2019}
Rahaman, N., Baratin, A., Arpit, D., Draxler, F., Lin, M., Hamprecht, F.A., Bengio, Y., Courville, A.: On the {Spectral} {Bias} of {Neural} {Networks} (2019).
\newblock ArXiv:1806.08734

\bibitem{raissi_physics-informed_2019}
Raissi, M., Perdikaris, P., Karniadakis, G.E.: Physics-informed neural networks: a deep learning framework for solving forward and inverse problems involving nonlinear partial differential equations.
\newblock J. Comput. Phys. \textbf{378}, 686--707 (2019)

\bibitem{shang_overlapping_2025}
Shang, Y., Heinlein, A., Mishra, S., Wang, F.: Overlapping {Schwarz} preconditioners for randomized neural networks with domain decomposition.
\newblock Comput. Methods Appl. Mech. Eng. \textbf{442}, 118011 (2025)

\bibitem{wu_comprehensive_2023}
Wu, C., Zhu, M., Tan, Q., Kartha, Y., Lu, L.: A comprehensive study of non-adaptive and residual-based adaptive sampling for physics-informed neural networks.
\newblock Computer Methods in Applied Mechanics and Engineering \textbf{403}, 115671 (2023)

\end{thebibliography}

\end{document}